\documentclass[12pt]{amsart} 
\usepackage{amssymb, epsfig} 
\usepackage{amsmath} 
\usepackage{latexsym} 
\numberwithin{equation}{section}  
\newcommand{\ol}{\overline}
\renewcommand{\a}{\alpha}  
\renewcommand{\b}{\beta}  
  
\newcommand{\G}{\Gamma}  
\renewcommand{\d}{\delta}

\renewcommand{\l}{\lambda}

\renewcommand{\o}{\omega}  
\renewcommand{\O}{\Omega}  
  
\newcommand{\s}{\sigma}  
  
\renewcommand{\t}{\tau}

\renewcommand{\P}{\mathcal P}  

\newcommand{\N}{{\mathbb N}}  
\newcommand{\R}{{\mathbb R}}

\newcommand{\A}{{\mathcal A}}  
\newcommand{\B}{{\mathcal B}}

\newcommand{\dis}{\displaysytle}

\newtheorem{theorem}{Theorem}[section]

\newtheorem{corollary}[theorem]{Corollary}  
  
\newtheorem{definition}[theorem]{Definition}  
\newtheorem{proposition}[theorem]{Proposition}  
\theoremstyle{definition}  
  
\theoremstyle{remark}  
\newtheorem{remark}[theorem]{Remark}

\title[Comparison principle for unbounded viscosity solutions]
{Comparison principle for unbounded viscosity solutions of degenerate
elliptic PDEs with gradient superlinear terms.}

\author[S. Koike and O. Ley]
{Shigeaki Koike \and Olivier Ley}

\address{ (S. Koike)  
Department of Mathematics, Saitama University\\
255 Shimo-Okubo, Sakura, Saitama 338-8570 Japan\\ {\tt 
skoike@rimath.saitama-u.ac.jp}}
\address{
(O. Ley) IRMAR, INSA de Rennes, F-35708 Rennes, France\\
{\tt olivier.ley@insa-rennes.fr}}

\thanks{S.K. was supported in part by Grant-in-Aid for Scientific Research 
(No. 20340026) of 
Japan Society for the Promotion of Science.
O.L. is partially supported by ANR BLANC07-3 187245, Hamilton-Jacobi and 
Weak KAM Theory. Parts of this work was done when S.K.  visited Tours 
University
and O.L.  visited Saitama University. The authors wish 
to acknowledge their hosts for hospitality.}



\begin{document}

\maketitle

\begin{abstract} 
We are concerned with fully nonlinear possibly degenerate elliptic partial 
differential equations (PDEs) with superlinear terms with respect to $Du$.
We prove several comparison principles among viscosity solutions 
which may be unbounded under some polynomial-type growth conditions. 
Our main result applies to PDEs with convex superlinear terms 
but we also obtain some results in nonconvex cases.
Applications to monotone systems of PDEs are given.
\end{abstract}

%%%%%%%%%%%%%%%%%%%%%%%%%%%%%%%%%%%%%%%%%%%%%%%%%%%%%%%%%%%%%%%%%%%%%%%%%%%%%

\def\l{\lambda}
\def\le{\left}
\def\ri{\right}
\def\la{\langle}
\def\ra{\rangle}
\def\R{\mathbb{R}}
\def\ti{\times}
\def\to{\rightarrow}
\def\fr{\frac}
\def\o{\omega}
\def\S{\mathcal{SSG}}
\def\SG{\mathcal{SG}}
\def\e{\varepsilon}
\def\M{\mathcal{M}}
\def\dis{\displaystyle}

\section{Introduction}

We are concerned with the comparison principle for viscosity solutions of 
 fully nonlinear elliptic partial differential equations:
\begin{equation}\label{1.1}
\l u+F(x,Du,D^2u)+H(x,Du)=f(x)\quad\mbox{in }\R^N,
\end{equation}
where $\l>0$, $F:\R^N\ti\R^N\ti S^N\to\R$, $H:\R^N\ti\R^N\to\R$ and 
$f:\R^N\to\R$ are 
given functions. 
Here $S^N$ denotes the set of $N\ti N$ symmetric matrices 
equipped with the standard order.

We will suppose that $F$ satisfies the standard hypothesis called structure 
condition. In particular, $F$ is degenerate elliptic, that is
\begin{equation}\label{degell}
F(x,\xi,X)\leq F(x,\xi,Y) \quad {\rm when} \ X\geq Y,
\quad x,p\in\R^N, \; X,Y\in S^N.
\end{equation}
On the contrary, we will suppose that the mapping $\xi\to H(x,\xi )$ has 
superlinear growth. 
A typical example is 
\begin{equation}\label{Htypical}
H(x,\xi)=\la A(x)\xi,\xi\ra^{q/2},
\end{equation}
where $q>1$, and $A:\R^N\to  S^N$.

When we consider unbounded solutions of
PDEs with superlinear growth terms in $Du$, 
we may not expect solutions to be unique in general. 
In fact, for $N=1$, the equation
\begin{equation}\label{1.2}
\l u-u''+|u'|^2=0\quad\mbox{in }\R
\end{equation}
admits at least two solutions; $u_1\equiv 0$ and 
$u_2(x)=-\fr{\l}{4}x^2-\fr{1}{2}$.

In \cite{alvarez96}, Alvarez introduced bounded-from-below 
 solutions to avoid $u_2$ in this case. 
He showed the uniqueness of strong bounded-from-below 
solutions of 
\begin{equation}\label{1.3}
u-\Delta u+|Du|^q=f(x)\quad \mbox{in }\R^N.
\end{equation}
We will mention this result after introducing some notations in 
Section~\ref{sec:2}.

We also refer to \cite{alvarez97} and \cite{hishii97} for comparison 
results, 
which yield the uniqueness among bounded-from-below viscosity solutions 
of Hamilton-Jacobi equations. 

On the other hand, 
 the uniqueness of unbounded viscosity 
solutions has been studied under certain growth condition on solutions. 
In this direction, H. Ishii \cite{hishii84} first 
established the comparison principle for 
unbounded viscosity solutions of Hamilton-Jacobi equations. 
For nonlinear elliptic PDEs, 
Aizawa-Tomita  \cite{AT1,AT2}, Crandall-Newcomb-Tomita \cite{CNT} and 
K. Ishii-Tomita \cite{KIT} 
obtained comparison results for 
unbounded viscosity solutions satisfying certain growth condition. 
However, unfortunately, we cannot apply these results to PDEs having 
variable 
coefficients to superlinear terms in $Du$. 
For instance, it seems difficult to treat typical $H$ as~\eqref{Htypical} 
unless $A$  is constant.

To avoid this technical difficulty, we will adapt a ``linearization'' 
technique, 
which Da Lio and the second author \cite{dll06} used to show the uniqueness 
of unbounded viscosity solutions of parabolic Bellman equations 
with quadratic nonlinearity.

More recently, we are informed that Barles and Porretta \cite{bp08} proved 
that \eqref{1.3} with $q=2$ 
admits at most one bounded-from-below solution if $f$ is bounded from 
below. 
In the case of \eqref{1.2}, $u_1$ is the only bounded-from-below solution.
However, their proof seems to be specific to \eqref{1.2} since
if we perturb this equation with a transport term as in
\begin{eqnarray} \label{hje3}
\l u-u''  +|u'|^2 +txu' =0\quad {\rm in} \ \R,
\end{eqnarray}
then there is at least two solutions $u_1\equiv 0$ and $u_2(x)
=-\fr{\l +2t}{4}x^2-\fr{\l +2t}{2\l}$. 
Thus, for $t<-\fr{\l}{2}$, $u_1$ and $u_2$ are 
bounded-from-below solutions of  (\ref{hje3}).
%This indicates that if $F$ contains a variable coefficient to $Du$, 
%then it seems hard to obtain comparison principle among bounded-from-below 
%viscosity solutions of (\ref{1.1}). 

In this paper, we study the comparison principle for viscosity 
solutions of (\ref{1.1}) under certain growth condition on $f$ 
and solutions.  We obtained two types of results depending
on whether $H(x,\xi)$ is convex  in $\xi$ or not.
The convex case is typically~\eqref{Htypical}  with positive 
$A(x)\in S^N$. Then we consider two nonconvex cases.
The first one is when $H(x,\xi)$  is convex in $\xi$ 
in some subset $\Omega_0\subset \Omega$ and
is concave in its complement. The second one
is when $H(x,\xi)$ is defined as a minimum of
convex Hamiltonians, that is,
$$
H(x,\xi)=\min\{ H_k(x,\xi) \ | \ k=1,\dots,m\},
$$ 
where $\xi\to H_k(x,\xi)$ is convex for $x\in\Omega.$
%But one of the main novelty is to able to deal with some cases for which $H(x,\xi)$ 
%changes its sign; few uniqueness results are known for PDEs with nonconvex $H$. (MAYBE REWRITE THIS SENTENCE). 
We will discuss a generalization of the above $H$, which appears in 
differential games (See \cite{nagai} for applications). 
Some applications to monotone systems of PDEs are also given.

Let us mention that we restrict ourselves to comparison principles
since it is the main ingredient to obtain existence and uniqueness
in the theory of viscosity solutions.

This paper is organized as follows: 
In Section~\ref{sec:2}, we give our hypothesis on $F$ and $H$. 
Section~\ref{sec:3} is devoted to the case when $H$ is strictly convex in 
$\xi$. 
We then discuss on the case when $H$ may be nonconvex in
Section~\ref{sec:4}. 
In section~\ref{sec:5}, we extend our results to  monotone systems.

%%%%%%%%%%%%%%%%%%%%%%%%%%%%%%%%%%%%%%

%%%%%%%%%%%%%%%%%%%%%%%%%%%%%%%%%%%%%%

\section{Preliminaries}

\label{sec:2}
First of all, we recall the definition of viscosity solutions of general 
PDEs:
\begin{equation}\label{general}
G(x,u,Du,D^2u)=0\quad\mbox{in }\R^N,
\end{equation}
where $G:\R^N\ti \R\ti \R^N\ti S^N\to \R$ is continuous. 
\begin{definition}
We call $u:\R^N\to \R$ a viscosity subsolution $($resp., supersolution$)$ 
of $(\ref{general})$ if 
for $\phi\in C^2(\R^N)$, 
$$
G(\hat x,u^*(\hat x),D\phi (\hat x),D^2\phi (\hat x))\leq 0$$
$$
(\mbox{resp., }
G(\hat x,u_*(\hat x),D\phi (\hat x),D^2\phi (\hat x))\geq 0)
$$
provided $u^*-\phi$ $($resp., $u_*-\phi)$ attains its local maximum 
$($resp., minimum$)$ at $\hat x\in \R^N$. 

We also call $u$ a viscosity solution of $(\ref{general})$ if it is both a 
viscosity sub- and supersolution of $(\ref{general})$. 
\end{definition}
Here $u^*$ and $u_*$ denote upper and lower semicontinuous envelopes of 
$u$, 
respectively.
We refer to \cite{cil92, barles94, bcd97, koike04} 
for their definitions, and the basic theory of 
viscosity solutions.

In order to explain our hypotheses below, we give a typical example:
\begin{equation}\label{example1}
u-\mbox{Tr}(\s (x)\s^T (x)D^2u)+\la b(x),Du\ra +\la 
A(x)Du,Du\ra^{\fr{q}{2}} 
=f(x)
\end{equation}
in $\R^N$, where $\s ,A:\R^N\to S^N$, and $b:\R^N\to \R^N$ are given 
functions. 
In this example, $G(x,\xi,X)=F(x,\xi,X)+H(x,\xi)-g(x)$ with
$F(x,\xi,X)=-\mbox{Tr}(\s (x)\s^T(x)X)+\la b(x),\xi\ra$, 
and $H(x,\xi)=\la A(x)\xi,\xi\ra^{q/2}$.

We denote by $\M$ the set of modulus of continuity; $m\in \M$ if 
$m(s)\to 0$ as $s\to 0^+$ and $m(s+t)\leq m(s)+m(t)$ for all $s,t>0.$

We present a list of hypothesis on $F$: 
The first one is a modification of the structure condition, under which 
we may consider (\ref{example1}) 
 when $\s$ and $b$ are locally Lipschitz continuous.
$$
\le\{\begin{array}{c}
\mbox{For }R>0, \mbox{ there exists }m_R\in \M\mbox{ such that}\\
F\le( x,\e^{-1}(x-y),X\ri)-
 F\le( y,\e^{-1}(x-y),Y\ri)\\
\leq m_R\le( 
|x-y|+\e^{-1}|x-y|^2\ri)\\
\mbox{provided }\e>0$, $x,y\in B_R\mbox{ and }(X,Y)\in S^N\ti S^N\mbox{ 
satisfies}\\
-\dis\fr{3}{\e}\le(
\begin{array}{cc}
I&O\\
O&I
\end{array}\ri)
\leq \le(
\begin{array}{cc}
X&O\\
O&-Y
\end{array}\ri)
\leq 
\dis\fr{3}{\e}
\le(
\begin{array}{cc}
I&-I\\
-I&I
\end{array}\ri)
.\end{array}
\ri.
\leqno{(F1)}$$
Here $B_r=\{ x\in\R^N \ | \ |x|<r\}$ and $B_r(x)=x+B_r$ for 
$r>0$ and $x\in\R^N$. 
Notice that (F1) implies the degenerate ellipticity (\ref{degell}).

We next suppose homogeneity of $F$ in $(\xi,X)\in \R^N\ti S^N$:
$$
F(x,\theta \xi,\theta X)=\theta F(x,\xi,X)\quad 
\mbox{for }\theta \geq 0,\,  x,\xi\in \R^N,\,  X\in S^N.\leqno{(F2)}
$$

To state further hypotheses, we introduce two subsets
of functions having superlinear growth %respectively strict sub-quadratic 
%and sub-quadratic
of order $r$;

A continuous function $h:\R^N\to \R$ belongs to $\S_r^\pm$
%(strict sub-quadratic growth)
if and only if
$$
\liminf_{|x|\to\infty}\fr{\pm h(x)}{|x|^r}\geq 0.
$$
Notice that $h\in \S_r^+$ (resp., $\S_r^-$)
if, for any $\e>0$, there exists $C_\e=C_\e(h)>0$ 
such that 
$$
h(x)\geq -\e |x|^r-C_\e\quad(\mbox{resp., } h(x)\leq \e |x|^r+C_\e)\quad 
\mbox{in }\R^N.
$$
We define $\S_r=\S_r^+\cap \S_r^-.$ Notice that
$h\in \S_r$ if and only if
$$
\lim_{|x|\to\infty}\frac{|h(x)|}{|x|^r}=0.
$$

A continuous function $h:\R^N\to \R$ belongs to $\SG_r^\pm$
%(sub-quadratic growth)
if and only if
$$
\liminf_{|x|\to\infty}\fr{\pm h(x)}{|x|^r}> -\infty.
$$
Notice that $h\in \SG_r^+$ (resp., $\SG_r^-$)
if, for any $\e>0$, there exist positive constants $\e=\e(h),C=C(h)$ 
such that 
$$
h(x)\geq -\e |x|^r-C\quad(\mbox{resp., } h(x)\leq \e |x|^r+C)\quad 
\mbox{in }\R^N.
$$
We define $\SG_r=\SG_r^+\cap \SG_r^-.$
Notice that, if a continuous function $h$ belongs to $\SG_r$, then
there exists $M>0$ such that, for all $x\in\R^N,$
$$
|h(x)|\leq M(1+|x|^r).
$$

The next assumptions indicate that the coefficients to the second and first 
derivatives are in $\S_2$ and $\S_1$, respectively.
$$
\le\{\begin{array}{c}
\mbox{There exists }
\s_0 :\R^N\to S^N\mbox{ such that } |\s_0|\in\S_1\mbox{ and }
\\
F(x,\xi,X)-F(x,\xi,Y)
\geq -\mbox{Tr}(\s_0(x)\s_0(x)^T(X-Y))\\
\mbox{ for }x,\xi\in\R^N, X,Y\in S^N.
\end{array}
\ri.
\leqno{(F3)}
$$
$$
\le\{\begin{array}{c}
\mbox{There exists }
b_0:\R^N\to\R\mbox{ such that  } |b_0|\in \S_1 \mbox{ and }
\\
|F(x,\xi,X)-F(x,\eta,X)|\leq b_0(x)|\xi-\eta|\\
\mbox{for }x,\xi,\eta \in\R^N, X\in S^N.
\end{array}
\ri.\leqno{(F4)}$$
We shall write $\P (x,X)=-\mbox{Tr}(\s_0(x)\s_0^T(x)X)$. 

We next give a list of hypothesis on $H:\R^N\ti\R^N\to\R$:
$$
\xi\in\R^N\to H(x,\xi)\mbox{ is convex for }x\in\R^N,
\leqno{(H1)}
$$
which will be violated in Section~\ref{sec:4} when we treat 
PDEs (\ref{example1}) with 
matrices $A(\cdot)$ which are not positive definite everywhere.  
Under $(H1)$, we need to suppose strict positivity and boundedness of 
$H$ with respect to $x\in\R^N$. For a fixed $q>1$,
$$
\le\{\begin{array}{c}
\mbox{There exist }\d\in C(\R^N) \mbox{ and }C_0>0\mbox{ such that }\d 
(x)>0,\\
\mbox{and }\d (x) |\xi|^q\leq H(x,\xi)\leq C_0 |\xi|^q
\mbox{ for }x,\xi\in\R^N.
\end{array}\ri.
\leqno{(H2)}$$
$$
H(x,\theta \xi)=\theta^q H(x,\xi)\quad\mbox{for }x,\xi\in\R^N,\,  \theta 
\geq 0.
\leqno{(H3)}
$$
We also suppose  continuity of $H$ in $x\in\R^N$. 
$$
\le\{\begin{array}{c}
\mbox{For }R>0,\mbox{ there exists }\o_R\in\M \mbox{ such that}\\
|H(x,\xi)-H(y,\xi)|\leq \o_R(|x-y|)|\xi|^q\\
\mbox{for }x,y\in B_R\mbox{ and }
\xi\in\R^N.
\end{array}
\ri.
\leqno{(H4)}
$$

In the sequel, we denote by $q'$ the conjugate of $q>1$;
$$
\frac{1}{q}+\frac{1}{q'}=1.
$$

Now, we shall come back to the result in \cite{alvarez96} for (\ref{1.3}). 
Roughly speaking, the comparison result in \cite{alvarez96} is as follows: 
if we suppose that
$$
f-g\in \S_{q'}\quad \mbox{for a nonnegative convex function }g:\R^N\to\R,$$
then the uniqueness holds among strong solutions in 
$W^{2,N}_{loc}(\R^N)\cap 
\S^+_{q'}$. 
Thus, if one restricts $f$ to be nonnegative and convex, 
then one does not need to suppose any 
growth condition on $f$ to obtain the comparison principle. 
In this paper, we generalize the uniqueness result
by assuming only that $f\in \S^+_{q'},$ i.e., $f$ may have
any growth from above and need not to be ``close'' to a
convex function. %(TO CHECK IF I AM RIGHT).

%In this paper, we only deal with the case of $g\equiv 0$ 
%in \cite{alvarez96}  but 
%we only assume $f\in \S^+_{q'}$ for some $t_f\geq 0$. 

%Although we could expect that the comparison principle holds 
%provided $f-g\in \S_{q'}^+$ for a nonnegative convex function $g$ 
%(more pricisely, for any $\e>0$ there is $C_\e>0$ such that 
%$f(x)\geq (1-\e)g(x)-\e |x|^{q'}-C_\e$), 
% we need to find strong supersolutions 
%of an extremal PDE derived from our linearization procedure. 
%We will come back this topics in a future work. 

%%%%%%%%%%%%%%%%%%%%%%%%%%%%%%%%%%%%%%%%%%%%%%%%%%%%%%%%%%%%%%%%%%%%
%%%%%%%%%%%%%%%%%%%%%%%%%%%%%%%%%%%%%%%%%%%%%%%%%%%%%%%%%%%%%%%%%%%%
%%%%%%%%%%%%%%%%%%%%%%%%%%%%%%%%%%%%%%%%%%%%%%%%%%%%%%%%%%%%%%%%%%%%

\section{Comparison principle}
\label{sec:3}
We  denote by $USC(\R^N)$ (resp. $LSC(\R^N)$) the set of upper (resp., 
lower) semicontinuous functions in $\R^N$. 
We first establish the comparison principle when given data are of $\SG_q$.

%%%%%%%%%%%%%%
\begin{theorem}\label{comparison1}Fix any $\l >0$. 
Assume that $(F1-4)$  and $(H1-4)$ hold. 
Let $u\in USC(\R^N)\cap \S^-_{q'}$ and $v\in LSC(\R^N)\cap \S_{q'}^+$ be, 
respectively, 
a viscosity subsolution and a viscosity supersolution of $(\ref{1.1})$. 
If $f\in \S_{q'}^+$, then for any $\l>0$, we have
$u\leq v$ in $\R^N$. 
\end{theorem}

%%%%%%%%%%%%

{\bf Proof.} 
{\it Step 1: Linearization procedure.} 

For $\mu\in (0,1)$, it is easy to verify that $u_\mu:=\mu u$ is a viscosity 
subsolution 
of
\begin{equation}\label{umu}
\l u_\mu +F(x,Du_\mu ,D^2u_\mu )+\mu^{1-q}H(x,Du_\mu)=\mu f 
(x)\quad\mbox{in }\R^N.
\end{equation}

We shall show that $w=w_\mu :=u_\mu -v$ is a viscosity subsolution of an 
extremal PDE
\begin{equation}\label{extremal}
\l w+\P(x,D^2w)-b_0(x)|Dw|
-\b_\mu |Dw|^q\leq (\mu -1)f(x)\quad\mbox{in }\R^N,
\end{equation}
where $\b_\mu:=(\fr{1-\mu}{2})^{1-q}C_0>0$.

For $\phi\in C^2(\R^N)$, we suppose that $w-\phi$ attains a local maximum 
at $\hat x\in\R^N$. 
We may suppose that $(w-\phi)(\hat x)=0>(w-\phi)(x)$ for $x\in B_r(\hat x)
\setminus \{ \hat x\}$ with a small $r\in (0,1)$.

Let $(x_\e, y_\e)\in B:=\overline{B}_r(\hat x)\ti \overline{B}_r(\hat x)$ 
be a maximum point of $u_\mu (x)-v(y)-(2\e)^{-1}|x-y|^2-\phi (y)$ over $B$. 
Since we may suppose $\lim_{\e\to 0}(x_\e,y_\e)=(\hat x,\hat x)$, and 
moreover 
$\lim_{\e\to 0}(u_\mu (x_\e),v(y_\e))=(u_\mu (\hat x),v(\hat x))$, 
it follows that $(x_\e,y_\e)\in {\rm int}(B)$ for small $\e.$ 
Hence, in view of Ishii's lemma (e.g. Theorem 3.2 in \cite{cil92}), 
setting $p_\e=\e^{-1}(x_\e-y_\e)$, 
we find $X_\e ,Y_\e\in S^N$ such that 
$(p_\e ,X_\e)\in \overline J^{2,+}u_\mu (x_\e)$, $(p_\e-D\phi 
(y_\e),Y_\e-D^2\phi (y_\e))
\in \overline J^{2,-}v(y_\e)$, and 
$$
-\fr{3}{\e}\le(\begin{array}{cc}
I&O\\
O&I
\end{array}\ri)
\leq \le(\begin{array}{cc}
X_\e &O\\
O&-Y_\e
\end{array}
\ri)
\leq 
\fr{3}{\e}\le(
\begin{array}{cc}
I&-I\\
-I&I
\end{array}\ri)
.$$
Thus, from the definition, we have
$$
\l u_\mu (x_\e)+F(x_\e ,p_\e,X_\e)+\mu^{1-q}H(x_\e,p_\e)\leq \mu f (x_\e)
$$
and
$$
\l v(y_\e)+F(y_\e,p_\e-D\phi (y_\e),Y_\e-D^2\phi (y_\e))+H(y_\e,p_\e-D\phi 
(y_\e))
\geq f(y_\e).$$
Since $(F3)$ and $(F4)$ imply 
$$
\begin{array}{rl}
&\P(y_\e,D^2\phi (y_\e))-b_0(y_\e)|D\phi (y_\e)|\\
\leq&F(y_\e,p_\e,Y_\e)-F(y_\e,p_\e -D\phi (y_\e),Y_\e-D^2\phi (y_\e)),
\end{array}
$$
by $(F1)$, 
we have 
$$
\begin{array}{rl}
&\l (u_\mu (x_\e)-v(y_\e))+\P(y_\e ,D^2\phi (y_\e))-b_0(y_\e)|D\phi 
(y_\e)|\\
\leq&
H(y_\e ,p_\e-D\phi (y_\e))-\mu^{1-q}H(x_\e,p_\e)+
\mu f (x_\e)-f(y_\e)\\
&+m_R(|x_\e-y_\e|+\e^{-1}|x_\e-y_\e|^2) 
,
\end{array}
$$
where $R=r+|\hat x|$.

We shall estimate the first two terms in the right hand side of the above. 
By $(H1)$, we have 
$$
H(y_\e,p_\e -D\phi (y_\e))\leq (\fr{1+\mu}{2})^{1-q}H(y_\e 
,p_\e)+(\fr{1-\mu}{2})^{1-q}
H(y_\e ,-D\phi (y_\e)).$$
Thus, due to $(H2)$ and $(H4)$,  we find $\o_R\in \M$ such that 
\begin{equation}\label{estimateH}
\begin{array}{rl}
&H(y_\e,p_\e-D\phi (y_\e))-\mu^{1-q}H(x_\e,p_\e)\\
\leq &-(\mu^{1-q}-(\fr{1+\mu}{2})^{1-q})\d (y_\e)|p_\e|^q +
\mu^{1-q}\o_R(|x_\e-y_\e|)|p_\e|^q\\
&+(\fr{1-\mu}{2})^{1-q}H(y_\e ,-D\phi (y_\e)).
\end{array}
\end{equation}
Since the positivity of $\d (\hat{x})$ implies 
 $\mu^{1-q}\o_R(|x_\e-y_\e|)\leq (\mu^{1-q}-(\fr{1-\mu}{2})^{1-q})
\d (y_\e)$  for small $\e>0$, we have 
$$\begin{array}{c}
\l (u_\mu (x_\e)-v(y_\e))+\P (y_\e ,D^2\phi (y_\e))-b_0(y_\e)|D\phi (y_\e)|
-{\b_\mu} |D\phi (y_\e)|^q\\
\leq \mu f(x_\e) -f(y_\e)+ m_R(|x_\e-y_\e|+\e^{-1}|x_\e-y_\e|^2),
\end{array}
$$
where ${\b_\mu}=(\fr{1-\mu}{2})^{1-q}C_0$. 
Therefore, sending $\e\to 0$ and using that $(2\e)^{-1}|x_\e-y_\e|^2\to 0,$ 
we have
$$
\l w(\hat x)+\P(\hat x,D^2\phi (\hat x))-b_0(\hat x)|D\phi (\hat x)|
-{\b_\mu} |D\phi (\hat x)|^q\leq (\mu -1)f(\hat x),$$
which proves that $w$ is a viscosity subsolution of~\eqref{extremal}.

{\it Step 2: Construction of smooth strict supersolutions of 
$(\ref{extremal}).$}

Let $\Phi (x)= (1-\mu)\{ C_1 + \a \la x\ra^{q'}\}$, where $\la x\ra 
=(1+|x|^2)^{1/2}$, 
and $C_1,  \alpha>0$ will be chosen later.

Note that 
$$D\la x\ra^{q'}=q'  \la x\ra^{q'-2}x,\mbox{ and }
D^2\la x\ra^{q'}=q'  \la x\ra^{q'-4}\left( \la x\ra^2 I+(q'-2)x\otimes 
x\right).$$
Since $\s_0,b_0\in \S_1$ and $f\in \S_{q'}^+$, for any $\e,\e'>0$, 
%setting $\Phi (x)=(1-\mu)( C_1+\a \la x\ra^{q'})$ for $\a>0$, 
we can find 
$C_\e=C_\e(\s_0,b_0)>0$ and $C_{\e'}=C_{\e'}(f)>0$  
(independent of $\a>0$) such that 
$$\begin{array}{rl}
&\P(x,D^2\Phi )-b_0(x)|D\Phi |+(1-\mu )f(x)\\
\geq& (1-\mu)\{-\a 
(\e \la x\ra^{q'}+C_\e\la x\ra^{q'-1})-\e' \la x\ra^{q'}-C_{\e'}\},
\end{array}
$$
and 
$$
-{\b}_\mu |D\Phi |^q\geq 
-(1-\mu)\a^q C_0'\la x\ra^{q(q'-1)}= 
%- (q')^q2^{q-1}C_0 (1-\mu)\a^q\la x\ra^{q(q'-1)}= 
-(1-\mu)\a^q C_0'\la x\ra^{q'},$$
where $C_0'=2^{q-1}(q')^q$. 
%(q')^q C_02^{q-1} 
Hence, we have
\begin{equation}\label{Phi}
\begin{array}{rl}
&\l\Phi +\P(x,D^2\Phi)-b_0(x)|D\Phi|-
{\b}_\mu |D\Phi|^q+(1-\mu)f(x)\\
\geq &(1-\mu)\{
\l C_1 +\a (\l-\e -C_\e\la x\ra^{-1}-\a^{q-1}C_0')\la x\ra^{q'}\\
&\hspace*{2cm} -\e'\la x\ra^{q'}-C_{\e'} \}.
\end{array}
\end{equation}
Fix $\e, \a\in (0,1)$ such that
$\e\leq \l/4$ and $\a^{q-1}C_0'\leq \l/4$. 
We then choose $\e'\leq \l\a /4$ to estimate 
the right hand side of the above from below by
$$
(1-\mu)\{\l C_1 -C_{\e'}
+\a (\fr{\l}{4}-C_\e\la x\ra^{-1})\la x\ra^{q'} \}.
$$
Hence, taking $C_1=\l^{-1}[ C_{\e'} +\max\{  C_\e \la x\ra^{q'-2}  \ |
\ \la x\ra\leq  4C_\e/\l\} ]+1$, 
we see that $\Phi$ satisfies 
\begin{equation}\label{Strict}
\l\Phi +\P(x,D^2\Phi)-b_0(x)|D\Phi|-{\b_\mu} |D\Phi|^q>
(\mu -1)f(x)\quad\mbox{in }\R^N.
\end{equation}

{\it Step 3: Conclusion.}

Since $w\in \S^-_{q'}$,  $w-\Phi$ takes its 
maximum at $\hat x\in \R^N$. 
Thus, we have
$$
\l w(\hat x)+\P(\hat x,D^2\Phi (\hat x))-b_0(\hat x)
|D\Phi (\hat x)|-{\b_\mu} |D\Phi (\hat x)|^q\leq (\mu-1)f(\hat x).$$

If $(w-\Phi )(\hat x)\geq 0$, then we get a contradiction to 
\eqref{Strict}. 
%because $\Phi$ is a classical strict supersolution of (\ref{extremal}). 
Hence,  we have
$$
w(x)\leq (1-\mu )(C_1+\a \la x\ra^{q'})\quad\mbox{for }x\in \R^N,$$
which concludes the assertion in the limit $\mu\nearrow 1$. 
\hfill\qed \\

Note that, if we suppose $\s_0\in \SG_1$ or 
$b_0\in \SG_1$ in (F3-4), then the 
comparison principle for \eqref{1.1} fails among solutions in $\SG_{q'}$ 
in general. 
In fact, we recall the example \eqref{hje3} stated in the Introduction.
In this example, $b_0\in \SG_1$ but does not belong to $\S_1$ 
unless $t=0$, and the
comparison obviously fails since one does not have uniqueness.

Also, if we consider 
\begin{equation}\label{Ex2}
u-(1+x^2)u''+|u'|^2=0\quad\mbox{in }\R,
\end{equation}
then it is easy to check that $v_1\equiv 0$ and 
$v_2(x)=\fr{1}{2}+\fr{1}{4}x^2$ 
are solutions of \eqref{Ex2} in $\SG_2$ but  $v_2\notin \S_2$. 
This nonuniqueness comes from $\s_0\in \SG_1$.

In \cite{KIT}, 
they may suppose that given functions belong to $\SG_1$
for the comparison principle. 
However, they need to suppose that $\l$ is large enough. 
We can extend their results
following the above arguments. 
$$
\le\{\begin{array}{c}
\mbox{There exists }\s_0:\R^N\to S^N\mbox{ such that }|\s_0|\in\SG_1\mbox{ 
and}\\
F(x,\xi,X)-F(x,\xi,Y)\geq -\mbox{Tr}(\s_0(x)\s_0(x)^T(X-Y))\\
\mbox{for }x,\xi\in\R^N, X,Y\in S^N,
\end{array}
\ri.
\leqno{(F3')}
$$
$$
\le\{\begin{array}{c}
\mbox{There exists }b_0:\R^N\to\R\mbox{ such that }|b_0|\in \SG_1\mbox{ 
and}\\
|F(x,\xi,X)-F(x,\eta,X)|\leq b_0(x)|\xi-\eta|\\
\mbox{for }x,\xi,\eta \in\R^N, X\in S^N.
\end{array}
\ri.\leqno{(F4')}$$

%%%%%%%%%%%%%
\begin{theorem}\label{comparison2}
Assume that $(F1,2)$, $(F3',4')$ and $(H1-4)$ hold. 
For $f\in \SG_{q'}^+$,  there exists $\l_0>0$ such that for $\l\geq\l_0$, 
if $u\in USC(\R^N)\cap \S^-_{q'}$ and $v\in LSC(\R^N)\cap \S_{q'}^+$ are, 
respectively, 
a viscosity subsolution and a viscosity supersolution of $(\ref{1.1})$, 
then $u\leq v$ in $\R^N$. 
\end{theorem}
%%%%%%%%%%%%%

{\bf Proof.} 
We do not need any change in Step 1 of proof of Theorem \ref{comparison1}.

In view of $(F3')$ and $(F4')$, we can get (\ref{Phi}) for some 
$\e,\e',C_\e,C_{\e'}>0$ which are not necessary small. 
Therefore, we can choose $\l_0>0$ such that for $\l\geq \l_0$, 
we can show $\Phi$ is a strict supersolution of (\ref{extremal}). 
The rest of proof can be done by the same argument.\hfill\qed
\smallskip

In the above Theorem, we need to assume that $u,-v\in \S^-_{q'}$ 
to be sure that $w-\Phi$ achieves a maximum in $\R^N$
(recall that $(1-\mu)$ in front of $\Phi$ is arbitrarily small).
If we are concerned with PDEs (1.1) without superlinear terms,
that is
\begin{equation}\label{3.4}
\l u+F(x,Du,D^2u)=f(x)\quad\mbox{in }\R^N,
\end{equation}
then we can obtain  slightly stronger results. 

%%%%%%%%%%%%%
\begin{proposition}\label{comparison3}
Assume that $(F1-4)$  holds. 
Let $u\in USC(\R^N)\cap \SG^-_{q'}$ and $v\in LSC(\R^N)\cap \SG_{q'}^+$ be, 
respectively, 
a viscosity subsolution and a viscosity supersolution of $(\ref{3.4})$. 
If $f\in \S_{q'}^+$, then 
$u\leq v$ in $\R^N$.
\end{proposition}
%%%%%%%%%%%%

{\bf Proof.} 
Following the argument in the proof of Theorem 3.1, we verify that 
$w:=u-v$ is a viscosity subsolution of 
\begin{equation}\label{3.5}
\l w+\P (x,D^2w)-b_0(x)|Dw|=0\quad\mbox{in }\R^N.
\end{equation}
Now, setting $\Phi (x)=\a\la x\ra^{q'}+C_1$ for $\a,C_1\geq 1$, we see that 
$\Phi$ satisfies 
\begin{equation}\label{3.6}
\l \Phi (x)+\P(x,\Phi (x))-b_0(x)|D\Phi (x)|\\
\geq (\l C_1-C_\e)+\a \la x\ra^{q'}(\l -\e \la x\ra^{-2}-\e\a^{-1}),
\end{equation}
where $\e>0$ is small enough so that the second term of the right hand side 
is positive. 
We then choose $C_1\geq C_\e/\l$ to show that $\Phi$ is a strict 
supersolution 
of (\ref{3.5}). 
Since we may take $\a$ large enough so that $w-\Phi$ attains its 
maximum at a point in $\R^N$, we conclude the proof.\hfill\qed
\smallskip

Finally, we treat the case when given functions are in $\SG_1$. 

\begin{proposition}\label{comparison4}
Assume that $(F1,2)$ and $(F3',4')$  hold. 
For $f\in \SG_{q'}^+$, there exists $\l_0>0$ such that 
if $u\in USC(\R^N)\cap \SG^-_{q'}$ and $v\in LSC(\R^N)\cap \SG_{q'}^+$ are, 
respectively, 
a viscosity subsolution and a viscosity supersolution of $(\ref{3.4})$, 
then $u\leq v$ in $\R^N$.
\end{proposition}

{\bf Proof.} 
As above, we can show (\ref{3.6}) but $\e>0$ may not be small. 
However, again, for large $\l>0$, 
we can show that $\Phi$ is a strict supersolution of (\ref{3.5}) 
when $\a,C_1$ are any large numbers. 
Thus, we can conclude the proof even for $w\in \SG^-_{q'}$.\hfill\qed

%%%%%%%%%%%%%%%%%%%%%%%%%%%%%%%%%%%%%%%%%%%%%%%%%%%%%%%%%%%%%%%%%%
%%%%%%%%%%%%%%%%%%%%%%%%%%%%%%%%%%%%%%%%%%%%%%%%%%%%%%%%%%%%%%%%%%%%
%%%%%%%%%%%%%%%%%%%%%%%%%%%%%%%%%%%%%%%%%%%%%%%%%%%%%%%%%%%%%%%%%%%%

\section{Non-convex $H$}
\label{sec:4}

In this section, we deal with some case when  $(H1)$ is not satisfied.

We denote by $\G\subset\R^N$ the zero-level set of $H(\cdot ,\xi)$ for all 
$\xi\in\R^N$;
$$
\G=\{ x\in \R^N \ | \ H(x,\xi)=0\mbox{ for any }\xi\in\R^N\}.$$
Our assumptions are as follows. For $\s_0$ in $(F3)$ and $b_0$ in $(F4)$, 
$$
\G\subset \{ x\in \R^N \ | \ \s_0(x)=0, \ b_0(x)=0\}.
\leqno{(A1)}
$$
Assumption $(A1)$ is a kind of degeneracy condition on the coefficients
of $F.$
$$
\le\{
\begin{array}{c}
\mbox{There exist open sets }\O^\pm\subset\R^N, \d^\pm\in C(\R^N)\mbox{ and}\\
C_0^\pm>0\mbox{ such that }\R^N =\G\cup \O^+\cup \O^-, 
\d^\pm (x)>0,\\
\d^\pm (x)|\xi|^q\leq \pm H(x,\xi)\leq C_0^\pm |\xi|^q
\mbox{ for }x\in\O^\pm, \xi\in\R^N,\\
\mbox{and }\xi\to \pm H(x,\xi)\mbox{ are convex for }x\in \O^\pm.\\
\end{array}
\ri.
\leqno{(A2)}
$$
It means that we can divide $\R^N\setminus \G$ into two open subsets:
$\O^+$ where $H(x,\cdot)$ is convex and $\O^-$ where $H(x,\cdot)$ is
concave.

When  $A(x)=a(x)I$ in (\ref{example1}) for some $a:\R^N\to\R$, 
$\O^\pm=\{ x\in \R^N \ | \ \pm a(x)>0\}$, and $\G=\{ x\in \R^N \ | \ 
a(x)=0\}$.  

We also suppose that $\s_0$ in $(F3)$ and $b_0$ in $(F4)$ satisfy that
$$
\s_0,b_0\in W^{1,\infty}_{loc}(\R^N).
\leqno{(A3)}
$$
Finally, we need some degeneracy condition for $H$ on $\Gamma.$
$$
\le\{\begin{array}{c}
\mbox{For each }x_0\in \G, \mbox{ there exist  }r,C_1>0\mbox{ such that}\\
|H(x,\xi )|\leq C_1|x-x_0|^q|\xi |^q\mbox{ for }x\in B_r(x_0).
\end{array}
\ri.
\leqno{(A4)}
$$

%%%%%%%%%%%%%%%%
\begin{theorem}\label{comparison3bis}
Assume that $(F1-4)$, $(H3,4)$ and $(A1-4)$ hold. 
Let $u\in USC(\R^N)\cap \S^-_{q'}$ and $v\in LSC(\R^N)\cap \S_{q'}^+$ be, 
respectively, 
a viscosity subsolution and a viscosity supersolution of $(\ref{1.1})$. 
If $f\in \S_{q'}^+$, then 
$u\leq v$ in $\R^N$. 
\end{theorem}
%%%%%%%%%%%%%%%%%%

{\bf Proof.} 
We first notice that the comparison principle holds if 
$\xi\to H(x,\xi)$ is concave instead of $(H1)$. 
In fact, we may take $w_\mu =u-\mu v$ for $\mu\in (0,1)$, and then 
we can follow the argument in the proof of Theorem \ref{comparison2}.

{\it Step 1: $u\leq f/\l \leq v$ on $\Gamma$.}

We only prove the first inequality since the second one can be shown  
similarly.
For $x_0\in \Gamma$, let $x_\e\in \overline B_1(x_0)$ be the 
maximum point of $u(x)-f(x_0)-(2\e)^{-1}|x-x_0|^2$ over 
$\overline{B}_1(x_0)$.
It is easy to see that $\lim_{\e\to 0}x_\e = x_0$; $x_\e\in B_1(x_0)$ for 
small 
$\e>0$.

It follows that we can write the viscosity inequality for the 
subsolution $u$ of \eqref{1.1} at $x_\varepsilon$ (see $e.g.$ 
\cite{cil92}): 
for any $\e>0$, there exists $X_\e\in S^{N}$
such that
\begin{equation} \label{ineq-mat732}
\left(p_\e,X_\e\right) \in\bar{J}^{2,+}u(x_\e), \quad  \mbox{with}
\quad  X_\e\leq \frac{3}{\e}I,
\end{equation}
where $p_\e=\e^{-1}(x_\e-x_0)$. 
We have
\begin{eqnarray*} 
\l u(x_\e)-\P(x_\e,X_\e)-b_0(x_\e)|p_\e|+H(x_\e,p_\e)\leq f(x_\e).
\end{eqnarray*}
By $(A3)$ and $(A4)$, we can find some constants 
$C_{\sigma, 1}, C_{b, 1}, C_1 >0$ such that, for $\e$ small enough, we
have
\begin{eqnarray*}
&& |\sigma_0 (x_\e)|\leq  C_{\sigma, 1}|x_\e -x_0|, \quad
|b_0 (x_\e)|\leq  C_{b, 1}|x_\e -x_0|, \\
&& {\rm and} \quad
|H(x_\e, p_\e)|\leq C_1 |x_\e -x_0|^q |p_\e|^q.
\end{eqnarray*}
It follows that there exists $C>0$ such that 
$$
\begin{array}{c}
\l u(x_\e)-C\le(\e^{-1}|x_\e-x_0|^2+\e^{-q}|x_\e-x_0|^{2q}\ri) \leq 
f(x_\e).
\end{array}
$$
Since $\lim_{\e\to 0}\e^{-1}|x_\e-x_0|^2=0$ and $\lim_{\e\to 
0}u(x_\e)=u(x_0)$, 
letting $\e\to 0$, we get
\begin{eqnarray*} 
\l u(x_0)\leq  f(x_0).
\end{eqnarray*}

{\it Step 2: Comparison on $\Omega^+\cup \G.$}

We can proceed exactly as in the convex case (Step 1 in the proof
of Theorem~\ref{comparison1}) 
to prove
that $w_\mu=\mu u-v$ (for $0<\mu <1$) is a subsolution of 
\eqref{extremal} in $\Omega^+$. 
Define $\Phi = (1-\mu )( C_1+\alpha \la x\ra^{q'})$ with the same
choice of constant $ \alpha, C_1$ as before. 
Notice that, with this choice, $\l \Phi\geq (\mu -1)f$ in $\R^N.$

Consider $\sup_{\Omega^+\cup\Gamma} ( w_\mu-\Phi_\mu )$. 
Since $w_\mu\in \S^-_{q'}$, 
this supremum is finite and is achieved at a point $\bar{x}$
which belongs to the closed set $\Omega^+\cup\G$. 
We distinguish two cases.

At first, if $\bar{x}\in \Omega^+$, then, arguing as in
the convex case (Step 2 in the proof of Theorem~\ref{comparison1}) we can 
write the 
viscosity inequality for $w_\mu$ using $\Phi$ as a test-function to 
show that the supremum is nonpositive.

Now, if $\bar{x}\in \Gamma,$
then, from Step 1, we get $u(\bar{x})\leq f(\bar{x})/\l \leq v(\bar{x})$
and therefore $w_\mu (\bar{x})\leq (\mu-1)f(\bar{x})/\l \leq
\Phi(\bar{x})$; thus the supremum is nonpositive. 
In both case, $w_\mu-\Phi\leq 0$. 
Letting $\mu\nearrow 1$,
we conclude $u\leq v$ in $\Omega^+\cup\Gamma$. 

{\it Step 3: Conclusion.}

To get the comparison in $\Omega^-\cup\Gamma,$ we use the fact that
we are in the concave case in $\Omega^-$. 
As noticed before, we can prove $u\leq v$ in $\O^-\cup \G$. 
\hfill\qed

%%%%%%%%%%%%%%%%%%%%%%%%%%%%%%%%%%%%%%%%%%%%%%%%%%%%%%%%%%%%%%%%%%%%%%%%%

%%%%%%%%%%%%%%%%%%%%%%%%%%%%%%%%%%%%%%%%%%%%%%%%%%%%%%%%%%%%%%%%%%%%%%%%%

In Introduction, we give a nonconvex $H:\Omega\ti\R^N\to \R$ defined by
\begin{eqnarray}\label{minconvex}
H(x,\xi)=\min\{ H_k(x,\xi ) \ | \ k=1,2,\ldots ,m\},
\end{eqnarray}
where $H_k$ is convex in $\xi$ and $m\in\N.$
We shall denote by $A$ the set $\{ 1,2,\ldots ,m\}$. 

\begin{theorem}\label{nonconvex2}
Assume that $(F1-4)$ holds, that $H$ in~\eqref{1.1}
is given by~\eqref{minconvex} and that
$(H1-4)$ holds for each $H_k$ with common 
$\d\in C(\R^N)$, $C_0>0$ and $\omega_R\in\M$ for $k\in A$. 
Let $u\in USC(\R^N)\cap \S^-_{q'}$ and $v\in LSC(\R^N)\cap \S^+_{q'}$ be, 
respectively, a viscosity subsolution and a viscosity supersolution 
of~\eqref{1.1}.
If $f\in \S^+_{q'}$, then for any $\l>0$, we have $u\leq v$ in $\R^N$. 
\end{theorem}

{\bf Proof. } 
It is enough to verify Step 1 in the proof of Theorem \ref{comparison1}. 
More precisely, we only need to check if (\ref{estimateH}) holds. 
We shall use the same notation in the proof of Theorem \ref{comparison1}. 
For any $\e>0$, we can choose $k_\e\in A$ such that 
$$
H(x_\e,p_\e)=H_{k_\e}(x_\e,p_\e).
$$
Hence, we have 
\begin{eqnarray*}
&& H(y_\e,p_\e-D\phi (y_\e))-\mu^{1-q}H(x_\e,p_\e)\\
&\leq &  H_{k_\e}(y_\e,p_\e-D\phi (y_\e))-\mu^{1-q}H_{k_\e}(x_\e,p_\e)\\
& \leq &-(\mu^{1-q}-(\fr{1+\mu}{2})^{1-q})\d 
(y_\e)|p_\e|^q+\mu^{1-q}\omega_R 
(|x_\e -y_\e |)|p_\e|^q\\
&&+(\fr{1-\mu}{2})^{1-q}C_0|D\phi (y_\e)|.
\end{eqnarray*}
Therefore, since the remaining proof is the same as in the proof of 
Theorem \ref{comparison1}, we conclude the proof. \hfill\qed
\medskip

We shall  generalize  the above $H$. 

Let $\A$ and $\B$ be compact metric spaces. 
For $\a\in \A$, $\b\in \B$, we consider continuous functions 
$\s ,\t :\R^N\ti \A\ti  \B\to M(N,n)$, where $M(N,n)$ denotes the set of 
$N\ti n$ real-valued matrices. 
For $\a\in \A$, $\b\in \B$, $a,b\in\R^n$, $x,\xi\in\R^N,$
we define $H^{\a,a}_{\b,b}:\R^N\ti\R^N\to 
\R$ by 
$$
H^{\a,a}_{\b,b}(x,\xi)=2\la \s (x,\a,\b)a-\t (x,\a,\b)b,\xi\ra -|a|^2+|b|^2.$$
We next set 
$$
\begin{array}{rcl}
H_{\b,b}(x,\xi )&=&\dis\sup_{\a\in \A,a\in\R^n}H^{\a,a}_{\b,b}(x,\xi)\\
&=&\dis\sup_{\a\in \A}\{|\s^T (x,\a,\b)\xi |^2-2\la \t (x,\a,\b)b,\xi\ra \} +|b|^2
\end{array}
$$
for $\b\in \B$, $b\in\R^n$ and $x,\xi\in \R^N.$ 
Finally, set
\begin{equation}\label{4.1}
\begin{array}{rcl}
H(x,\xi)&=&\dis\inf_{\b\in \B,b\in\R^n}H_{\b,b}(x,\xi)\\
&=&\dis\inf_{\b\in \B}\dis\sup_{\a\in \A}\{|\s^T (x,\a,\b)\xi|^2-
|\t^T (x,\a,\b)\xi|^2\}.
\end{array}
\end{equation}
%%%
Defining $S(x,\a,\b)=\s (x,\a,\b) \s^T (x,\a,\b),$ 
$T(x,\a,\b)=\t (x,\a,\b) \t^T (x,\a,\b)\in~S^N,$ for 
$x\in \R^N$ and $(\a,\b)\in\A\ti \B$, we give a condition on $S, T$ so that $H$ satisfies $(H2)$. 

$$
\le\{\mbox{
\begin{tabular}{c}
There are $\d\in C(\R^N)$ and $C_0>0$ such that\\
\begin{tabular}{ll}
(i)&$\d (x)>0$ for $x\in\R^N$,\\
(ii)&for any $x\in\R^N$ and $\b\in \B$, there exists $\a_{\b,x}\in \A$\\
&satisfying $S(x,\a_{\b,x},\b)- T(x,\a_{\b,x},\b)\geq \d (x) I$,\\
(iii)&for any $x\in\R^N$, there exists $\b_x\in\B$ satisfying\\
&$\sup_{\a\in\A}| S(x,\a,\b_x)|\leq C_0$.
\end{tabular}
\end{tabular}}\ri.
\leqno{(H2')}$$
Assuming that $S, T:\R^N\ti \A\ti \B\to S^N$ satisfy $(H2')$, we 
easily verify that the above $H$ satisfies $(H2)$ and $(H3)$ with $q=2$. 
In fact, for $x,\xi\in\R^N$, we choose $\b_x=\b_{x,\xi}\in \B$ such that 
$H(x,\xi)=\sup_{\a\in A}\{ |\s^T (x,\a,\b_x)\xi|^2-|\t^T (x,\a,\b_x)\xi|^2\}$. 
Thus, by $(H2')$, we can find $\a_x=\a_{x,\xi}\in \A$ such that 
$$\begin{array}{rcl}
H(x,\xi)&\geq &|\s^T (x,\a_x,\b_x)\xi|^2-|\t^T (x,\a_x,\b_x)\xi|^2\\
&=&\la (S(x,\a_x,\b_x)-T(x,\a_x,\b_x))\xi,\xi\ra \\
&\geq& \d (x)|\xi|^2.
\end{array}
$$
The other inequality is trivial by (iii) of $(H2')$. 
Furthermore, assuming that 
$$\le\{\begin{array}{ll}&\mbox{for }R>0, \mbox{ there are }C_R>0\mbox{ and }\hat\o_R\in\M\mbox{ such that}\\
\text{(i)}&| \s(x,\a,\b)|+|\t (x,\a,\b)|\leq C_R\mbox{ for }x\in B_R\mbox{ and}\\
&(\a,\b)\in\A\ti\B ,\\
\text{(ii)}&|\s (x,\a,\b)-\s (y,\a,\b)| +|\t (x,\a,\b)-\t (y,\a,\b)|\\
&\leq \hat\o_R(|x-y|) \mbox{ for }x,y\in B_R$ and $(\a,\b)\in \A\ti\B ,
\end{array}\ri.
\leqno{(H4')}$$
we can show that $(H4)$ holds with some $\o_R\in\M$. 

Now, we can state the comparison principle for the above $H$ in \eqref{1.1}. 
Since we can prove it with the same argument as in the proof of Theorem \ref{nonconvex2}, 
we leave it to the readers. 

\begin{corollary}
Assume that $(F1-4)$ holds, that $H$ in~\eqref{1.1}
is given by~\eqref{4.1} and that
$(H2'), (H4')$ hold. 
Let $u\in USC(\R^N)\cap \S^-_{2}$ and $v\in LSC(\R^N)\cap \S^+_{2}$ be, 
respectively, a viscosity subsolution and a viscosity supersolution 
of~\eqref{1.1}.
If $f\in \S^+_{2}$, then for any $\l>0$, we have $u\leq v$ in $\R^N$. 
\end{corollary}

In particular, we shall suppose that $\s$ and $\t$ are, respectively, 
independent of $\a$ and $\b$. 
Then, it is easy to see 
$$
H(x,\xi)=\min_{\b\in \B}|\s^T (x,\b)\xi|^2-\min_{\a\in \A}|\t^T (x,\a)\xi|^2.
$$
Since it is straightforward to restate the hypothesis $(H2')$ and $(H4')$ in this case, 
we leave it to the readers. 
 
\begin{remark}
We may give some generalizations of
Theorems~\ref{comparison3bis} and \ref{nonconvex2} 
to PDEs with coefficients in $\SG$ instead
of $\S$ as it was done at the end of Section~\ref{sec:3}. 
\end{remark}

%%%%%%%%%%%%%%%%%%%%%%%%%%%%%%%%%%%%%%%%%%%%%%%%%%%%%%%%%%%%%%%%%%%%%%%%%

\section{Monotone systems}
\label{sec:5}

In this section, we establish the comparison principle to monotone systems 
of elliptic PDEs, which were introduced in \cite{IK91}. 

For a given integer $m\geq 2$, we set $A=\{ 1,2,\ldots,m\}$. 
We consider systems of PDEs: for $k\in A$, 
\begin{equation}\label{monosys}
F_k(x,u,Du_k,D^2u_k)+H_k(x,Du_k)=f_k(x)\quad\mbox{in }\R^N,
\end{equation}
where $u=(u_1,u_2,\ldots ,u_m):\R^N\to\R^m$ is an unknown function, and 
$F_k:\R^N\ti \R^m\ti \R^N\ti S^N\to \R$, $H_k:\R^N\ti\R^N\to\R$, 
$f_k:\R^N\to\R$ ($k\in A$) are given functions.

First of all, we recall the definition of viscosity solutions of general 
systems of PDEs: for $k\in A$, 
\begin{equation}\label{generalsys}
G_k(x,u,Du_k,D^2u_k)=0\quad\mbox{in }\R^N,
\end{equation}
where $G_k:\R^N\ti \R^m\ti \R\ti \R^N\ti S^N\to \R$ is continuous.

%%%%

\begin{definition}
We call $u=(u_k):\R^N\to \R^m$ a viscosity subsolution $($resp., 
supersolution$)$ 
of $(\ref{generalsys})$ if 
for $\phi\in C^2(\R^N)$ and $k\in A$, 
$$
G_k(\hat x,u^*(\hat x),D\phi (\hat x),D^2\phi (\hat x))\leq 0$$
$$
(\mbox{resp., }
G_k(\hat x,u_*(\hat x),D\phi (\hat x),D^2\phi (\hat x))\geq 0)
$$
provided $(u_k)^*-\phi$ $($resp., $(u_k)_*-\phi)$ attains its local maximum 
$($resp., minimum$)$ at $\hat x\in \R^N$. 

We also call $u$ a viscosity solution of $(\ref{generalsys})$ if it is both 
a 
viscosity sub- and supersolution of $(\ref{generalsys})$. 
\end{definition}

%%%%%

We will suppose that $F:=(F_1,F_2,\ldots ,F_m)$ is monotone as in 
\cite{IK91}: 
$$
\le\{\begin{array}{c}
\mbox{There exists }\l>0\mbox{ such that}\\
\mbox{if } r=(r_k), s=(s_k)\in\R^m, \  (x,\xi, X)\in \R^N\ti
\R^N\times S^N \mbox{ and }\\
 \max_{k\in A}(r_k-s_k)= r_j-s_j\geq 0
\mbox{ for } j=j(r,s,x,\xi ,X)\in A,\\
\mbox{then } F_j(x,r,\xi ,X)-F_j(x,s,\xi ,X)\geq \l (r_j-s_j).
\end{array}\ri.
\leqno{(M)}
$$
We will suppose that every $F_k=F_k(x,r,\xi,X)$ in $F=(F_k)$
satisfies $(F1)$ with a modulus $m_{R,k}$
uniformly for $|r|\leq R$; moreover it satisfies 
$(F3)$ and $(F4)$ with some 
$\sigma_k\in \S_1$ and $b_k\in \S_1,$ respectively.
Assumption $(F2)$ is replaced with
$$
F(x,\theta r,\theta \xi,\theta X)=\theta F(x,r,\xi,X)\quad 
\mbox{for }\theta \geq 0,  x,\xi\in \R^N,  r\in\R^m,
X\in S^N.\leqno{(F2')}
$$
We set
$\P_k (x,X)=-\mbox{Tr}(\s_k(x)\s_k^T(x)X)$. 
In the same way, we will assume that
$H_k$ satifies $(H1)$--$(H4)$ with common $\d\in C(\R^N),$ 
$q>1,$ and $\o_R$ (though we 
may allow them to depend on $k\in A$).

%%%%%%%%

\begin{theorem}\label{comparison1bis}
Assume that $(M)$, $(F1,2',3,4)$ hold for $F_k$ 
and $(H1-4)$ hold for $H_k$ $(k\in A)$.

Let $u_k\in USC(\R^N)\cap \S^-_{q'}$ and $v_k\in LSC(\R^N)\cap \S_{q'}^+$, 
$u=(u_k)$ and $v=(v_k)$ be, respectively, 
a viscosity subsolution and a viscosity supersolution of $(\ref{monosys})$. 
If $f_k\in \S_{q'}^+$ for $k\in A$, then 
$u_k\leq v_k$ in $\R^N$ for $k\in A$. 
\end{theorem}

%%%%%%%%

{\bf Proof.} 
First of all, by $(F2)$, $(H1)$ and $(H3)$, 
we verify that $u_\mu=(u_{\mu, k})=(\mu u_k)$ ($\mu \in (0,1)$) 
is a viscosity subsolution of 
$$
F_k(x,u_\mu ,Du_{\mu, k},D^2u_{\mu, k})+\mu^{1-q}H_k(x,Du_{\mu, k})\leq 
\mu f_k(x)\quad\mbox{in }\R^N.$$

{\it Step 1: Linearization.}
Set $w(x)=\max_{k\in A}(u_{\mu,k}-v_k)(x)$ for $x\in \R^N$. 
We shall verify that $w$ is a viscosity subsolution of
\begin{eqnarray*}\label{extremal-system}
\l w+\min_{k\in A}\{ \P_k(x,D^2w)-b_k(x)|Dw|-\b_\mu |Dw|^q-(\mu 
-1)f_k(x)\}=0
\end{eqnarray*}
in $\R^N,$
where $\b_\mu =\left(\frac{1-\mu}{2}\right)^{1-q} C_0.$
We argue as in the proof of Theorem \ref{comparison1}
assuming that, for  a fixed $\phi\in C^2(\R^N)$,
$w-\phi$ attains a strict local maximum at $\hat x\in\R^N.$
Setting $B:=\ol B_r(\hat x)\ti \ol B_r(\hat x),$
up to extract subsequences, we can suppose that
\begin{eqnarray}
\label{bla14} 
&& \dis\max_{x,y\in B}\dis\max_{k\in A}
 \{ u_{\mu, k}(x)-v_k(y)-(2\e)^{-1}|x-y|^2-\phi (y)\}\\\nonumber
&=& u_{\mu, j(\e)}(x_\e)-v_{j(\e)}(y_\e)-(2\e)^{-1}|x_\e-y_\e|^2-\phi(y_\e)
\end{eqnarray}
$x_\e, y_\e\to \hat x$ and $u_{\mu, j(\e)}(x_\e)-v_{j(\e)}(y_\e)\to w(\hat 
x).$
Moreover, since the set $A$ is finite, we may suppose that $j(\e)=j$
is independent of $\e.$

As in the proof of Theorem \ref{comparison1}, since there are $X_{j,\e}, 
Y_{j,\e}\in 
S^N$ such that $(p_\e,X_{j,\e})\in \ol J^{2,+}u_j(x_\e)$, $(p_\e -D\phi 
(y_\e),
Y_{j,\e}-D^2\phi (y_\e))\in \ol J^{2,-}v_j(y_\e)$, and 
the matrix inequalities in $(F1)$ hold with $(X_{j,\e},Y_{j,\e}),$ 
we have 
\begin{eqnarray}\label{form5.4}
&& F_j(y_\e,u_\mu (x_\e),p_\e,Y_{j,\e}) \\
&\leq&  F_j(x_\e,u_\mu (x_\e),p_\e,X_{j,\e})+  
m_R(|x_\e-y_\e|+\e^{-1}|x_\e-y_\e|^2),
\nonumber
\end{eqnarray}
where $R=r+|\hat x|$. 

Moreover, by $(F3)$ and $(F4)$, we have 
\begin{eqnarray}\label{form5.5}
&& F_j(y_\e,v(y_\e),p_\e-D\phi (y_\e),Y_{j,\e}-D^2\phi (y_\e))\\
&\leq&  F_j(y_\e,v(y_\e),p_\e,Y_{j,\e})-\P_j(y_\e,D^2\phi (y_\e))
+b_j(y_\e)|D\phi (y_\e)|.
\nonumber
\end{eqnarray}
From \eqref{bla14}, we note that
\begin{eqnarray*}
\max_{k\in A} (u_{\mu, k}(x_\e)-v_{k}(y_\e)) &=& u_{\mu, 
j}(x_\e)-v_{j}(y_\e)
\end{eqnarray*}
and therefore, by $(M),$ we have
\begin{eqnarray}\label{form5.10}
&&\l (u_{\mu, j}(x_\e)-v_{j}(y_\e))\\\nonumber
&\leq & F_j(y_\e,u_\mu (x_\e),p_\e,Y_{j,\e})-  
F_j(y_\e,v(y_\e),p_\e,Y_{j,\e}).
\end{eqnarray}

On the other hand, from the definition, we have 
$$
F_j(x_\e,u_\mu (x_\e),p_\e,X_{j,\e})+\mu^{1-q}H(x_\e,p_\e)\leq \mu 
f_j(x_\e),
$$
and
$$
F_j(y_\e,v(y_\e),p_\e-D\phi (y_\e),Y_{j,\e}-D^2\phi 
(y_\e))+H(y_\e,p_\e)\geq f_j(y_\e).
$$
Thus, following the same calculations for $H_j$ as in 
Theorem~\ref{comparison1}, 
by \eqref{form5.4}, \eqref{form5.5} and \eqref{form5.10}, we have
\begin{eqnarray}\label{form5.6}
&& \l (u_{\mu, j}(x_\e)-v_{j}(y_\e)) \\\nonumber
&& +\P_{j}(y_\e,D^2\phi (y_\e))- b_{j}(y_\e)|D\phi (y_\e)|
   -\b_\mu |D\phi (y_\e)|^q \\\nonumber
&& - \mu f_{j}(x_\e)+f_{j}(y_\e)\\\nonumber
&\leq& m_R (|x_\e-y_\e|+\e^{-1}|x_\e-y_\e|^2)
\end{eqnarray}
for small enough $\e>0$. 
Hence, sending $\e\to 0$ in \eqref{form5.6}, we obtain
the desired extremal PDE
$$
\l w(\hat x)+\min_{k\in A}\{ \P_k(\hat x,D^2\phi (\hat x))-b_k(\hat x)
|D\phi (\hat x)|-\b_\mu |D\phi (\hat x)|^q-(\mu -1)f_k(\hat x)\}\leq 0.
$$

{\it Step 2: Conclusion.}
Consider the same function $\Phi$ from the proof of
Theorem~\ref{comparison1}. 
We can choose the constant $\alpha ,C_0>0$ in order that
$\Phi$ is a strict supersolution of  the previous extremal
PDE. The conclusion follows.
\hfill\qed

\begin{remark}
As in the previous sections, we may give some generalizations of
Theorem~\ref{comparison3bis} to PDEs with coefficients in $\SG$ instead
of $\S$ and for nonconvex Hamiltonians $H_k$ satisfying assumptions
like (A1)--(A4) on some subsets $\O_k^\pm,$ $\Gamma_k.$ The proof
combines techniques developed in Section~\ref{sec:3} and~\ref{sec:4}, so 
we skip it. 
\end{remark}

%%%%%%%%%%%%%%%%%%%%%%%%%%%%%%%%%%%%%%%%%%%%%

%%% BIBLIO

%%%%%%%%%%%%%%%%%%%%%%%%%%%%%%%%%%%%%%%%%%%%%

%%%% AVEC BIBTEX

%\bibliographystyle{plain} 

%\bibliography{biblio} 

%\end{document}

%%%%%%%%%%%%%%%%%%%%%%%%%%%%%%

\end{document}